\documentclass[a4paper]{amsart}
\usepackage{graphicx,amssymb,amsmath}
\usepackage{hyperref}
\def\eps{\varepsilon}
\def\R{\mathbb R}
\def\Z{\mathbb Z}
\def\T{{\mathbb T^{1}}}

\newtheorem{theorem}{Theorem}
\newtheorem{remark}{Remark}
\newtheorem{proposition}{Proposition}
\newtheorem{corollary}{Corollary}
\newtheorem{lemma}{Lemma}

\begin{document}
\title{Bernoulli problem for rough domains}
\author[F. Bouchon and L. Chupin]{Fran\c{c}ois Bouchon and Laurent Chupin}
\address[F. Bouchon and L. Chupin]{Clermont Universit\'e, Universit\'e Blaise-Pascal, Laboratoire de Math\'ematiques, BP 10448, F-63000 CLERMONT-FERRAND FRANCE}
\address[F. Bouchon and L. Chupin]{CNRS, UMR 6620, LM, F-63171 AUBI\`ERE FRANCE}
\email{Francois.Bouchon@math.univ-bpclermont.fr}
\email{Laurent.Chupin@math.univ-bpclermont.fr}

\begin{abstract}
We consider the exterior free boundary Bernoulli problem in the case of a rough given domain. 
An asymptotic analysis shows that the solution of the initial problem can be approximated by the solution of a non-rough Bernoulli problem at order 2.
Numerical tests confirm these theoretical results.
\end{abstract}
\maketitle

\section{Introduction}

The free Boundary Bernoulli problem is a model problem with many applications in engineering sciences such as fluid mechanics (see \cite{F}) or electromagnetics (see \cite{C}, \cite{D} and references therein). 
It consists in an overdetermined elliptic problem on a free domain whose solution is the domain itself as well as the potential (see \cite{FR} for a review). More precisely, let~$\Omega$ be a bounded open set, we look for a domain~$A$ containing~$\overline\Omega$ and a function~$u$ defined on $A\setminus\overline\Omega$ such that:
\begin{equation}\label{pbintro}
\left\{
\begin{array}{rcll}
\Delta u&=&0&\mbox{ in }A\setminus\overline\Omega, \\
u&=&1&\mbox{ on }\partial\Omega, \\
u&=&0&\mbox{ on }\partial A,\\
|\nabla u|&=&\lambda&\mbox{ on }\partial A.
\end{array}
\right.
\end{equation}
This problem is known as ``exterior Bernoulli problem'', since the boundary of the solution domain $\partial A$ is exterior to the given domain $\Omega$. The interior Bernoulli problem is similar, with $A\subset\Omega$ and then $u$ defined (and harmonic) on $\Omega\setminus\overline A$.

Theoretical questions have been addressed concerning existence or uniqueness of solution as well as geometric properties.
In \cite{B}, Beurling proposed a classification of solutions (elliptic, parabolic and hyperbolic) and introduced the method of sub- and super-solution to prove the existence and uniqueness of the elliptic solution in two-dimensions. This method was later adapted by Acker in \cite{A} to get the same result in higher dimension in the case of a convex given domain $\Omega$. The same result was also proved by Alt and Caffarelli in \cite{AC} using variational methods.
Henrot and Shahgholian showed in \cite{HSh} that, if $\Omega$ is convex, then the solution domain $A$ is convex. They also give a positive answer to this question for the $p$-Laplacian case in~\cite{HSh2} and~\cite{HSh3}.

Acker and Meyer considered the more general starlike case in~\cite{AM}, and showed that the solution domain is starlike and elliptic. Using both geometric arguments and variational methods, they also proved that the solution domain depends continuously on the data (see also~\cite{HHS}).

Another class of works consider the numerical solution of this free boundary problem, most of them consisting in iterative algorithms which build a sequence of domains converging to the solution domain. A fixed point type approach have been introduced by Flucher and Rumpf in~\cite{FR}, and later combined with a levelset approach in~\cite{BCT}. 
Some numerical methods based on shape functional minimization have also been developped in~\cite{BBPST}, \cite{HKKP} and~\cite{IKP}, and for the three dimensional case in~\cite{H}.\\

In most situations of physical relevance, the given domain~$\Omega$ is issued from industrial situations and its boundary is generally rough.
Such geometry can not be described in detail: either the precise shape of the roughness is unknown, or its spatial variations are too small for computational grids. Therefore, one may only hope to account for the averaged effect of the irregularities.

Such an approach is usually used for Stokes flows: the irregular boundary is replaced by an artificial smoothed one, and an artificial boundary condition (a wall law) is prescribed there, that should reflect the mean impact of the roughness, see for instance~\cite{APV98,ABLS01,ChupinMartin12,JM03,NNRM}.

Roughly speaking, to obtain a "wall law" the method consists in correcting the flow without roughnesses using a contribution located in the neighborhood of the actual boundary.
The mathematical proofs are then based on the justification of the "local" correction.\\

The goal of the present article is to obtain a kind of "wall law" for a Bernoulli problem when the given domain~$\Omega$ is assumed to be rough.
For the Bernoulli problem the main difficulty is that the unknown is not only a function but also the domain in which this function is defined.
We then propose to study a simple case: we consider that the rough domain is a perturbation of a ball in~$\R^2$, and that the perturbation is given as a periodic function of amplitude and period equal to~$\eps$.
We discuss in the conclusion, see section~\ref{concl} some possible extensions of these results.\\

The rest of this paper is composed of three main parts.
Section~\ref{stat} is devoted to the presentation of the problem and its mathematical framework. We also give the main result.
The proof of this result is detailed in Section~\ref{proof} whereas some numerical simulations are given in Section~\ref{numerics}.
Finally, some technical results are postponed in Appendix.

\section{Statement of the problem and main results}\label{stat}

The Bernoulli problem that we consider reads: 
given a constant $\lambda>0$ and a domain~$\Omega_\eps\subset \R^2$, find a domain~$A_\eps$ containing~$\Omega_\eps$ and function~$u_\eps$ defined on~$A_\eps \setminus \overline{\Omega_\eps}$ such that:
\begin{equation}\label{pbberpert}
\left\{
\begin{array}{rcll}
\Delta u_\eps&=&0&\mbox{ in }A_\eps\setminus \overline{\Omega_\eps}, \\
u_\eps&=&1&\mbox{ on }\partial\Omega_\eps, \\
u_\eps&=&0&\mbox{ on }\partial A_\eps,\\
|\nabla u_\eps|&=&\lambda&\mbox{ on }\partial A_\eps.
\end{array}
\right.
\end{equation}
The parameter $\eps>0$ models the size of the roughness of the domain~$\Omega_\eps$. More precisely we consider a domain~$\Omega_\eps$ defined using polar coordinates:
\begin{equation}
\Omega_\eps = \Big\{ (r,\theta)\in \R_+\times \T ~;~ r < 1 - \eps h\Big( \frac{\theta}{\eps} \Big) \Big\},
\end{equation}
where $\T=\R/(2\pi \Z)$ denotes the torus, and where $h$ is a Lipschitz function on~$\T$.
\par
Note that a function $f$ defined on $\T$ can be viewed as a $2\pi$-periodic function $\widetilde f$ defined on $\R$.
Hence the parameter~$\eps$ must be the inverse of an integer number: $\frac{1}{\eps}\in \mathbb N$, and it must be small enough: $\frac{1}{\eps}>\displaystyle \max_\T h$.
We define the measure $\mathrm d\Theta$ as the measure quotient, that is such that
$$\int_\T f(\Theta) \, \mathrm d\Theta = \int_0^{2\pi} \widetilde f(t) \frac{\mathrm dt}{2\pi}.$$

We note $\Omega_0$ the domain without roughness, that is $\Omega_0=\mathcal B(0,1)$.
The function~$h$ is assumed to be a lipschitz and nonnegative function. That implies $\Omega_\varepsilon\subset\Omega_0$ for all $\varepsilon>0$.
Note that using these polar coordinates~$(r,\theta)$ the Laplacian operator reads
$$\Delta = \frac{\partial^2}{\partial r^2} + \frac{1}{r} \frac{\partial}{\partial r} + \frac{1}{r^2} \frac{\partial^2}{\partial \theta^2}.$$
For such a context, it is known that the Bernoulli problem~\eqref{pbberpert} has a unique solution, see~\cite{AM}. In fact this existence results from the fact that the domain~$\Omega_\eps$ is starlike with respect to all points in a ball $\mathcal B(0,\delta)$.\\

In order to approach this problem~\eqref{pbberpert} for $\eps \to 0$ using a Bernoulli problem on the non-rough domain~$\Omega_0$, the most natural idea is to introduce the following Bernoulli problem: find a domain~$A_0$ containing~$\Omega_0$ and function~$u_0$ defined on~$A_0 \setminus \overline{\Omega_0}$ such that
\begin{equation}\label{pbber0}
\left\{
\begin{array}{rcll}
\Delta u_0&=&0&\mbox{ in }A_0\setminus \overline{\Omega_0}, \\
u_0&=&1&\mbox{ on }\partial\Omega_0, \\
u_0&=&0&\mbox{ on }\partial A_0,\\
|\nabla u_0|&=&\lambda&\mbox{ on }\partial A_0.
\end{array}
\right.
\end{equation}
This approach is mathematically justified by the continuity results of~\cite{AM} (see also~\cite{HHS}). More precisely we have
\begin{proposition}
The solution~$(A_0,u_0)$ of the Bernoulli problem~\eqref{pbber0} approaches the solution~$(A_\eps,u_\eps)$ of the Bernoulli problem~\eqref{pbberpert} at order one:
\begin{equation}
\label{thm0}
\mathcal D_{\mathcal H}(A_0,A_\eps)=\mathcal O(\eps).
\end{equation}
\end{proposition}
The proof of this proposition directly follows from Theorem 3.9 in \cite{AM}.
More precisely it expresses that
$d_2(\partial A_0,\partial A_\varepsilon)\le d_2(\partial \Omega_0,\partial \Omega_\varepsilon)$
where $d_2$ is the metric defined by:
$d_2(\Gamma_1,\Gamma_2)=\sup\{|\ln(\gamma)|, \gamma\in\R_+ ~;~ \gamma\Gamma_1\cap\Gamma_2\neq\emptyset\}$.
For $0<\delta<M$, we introduce the set
\begin{align*}
\mathcal K=\{ \partial\Omega ~\text{where}~ \Omega\subset\mathcal B(0,M)& ~\text{is a starlike domain  with respect to all points in }\mathcal B(0,\delta)\}.
\end{align*}
We observe that there exists a choice of~$\delta$, $M$ and~$\eps_0$ such that $\partial \Omega_\varepsilon\in\mathcal K$ and $\partial A_\varepsilon\in\mathcal K$ for all $0\le\varepsilon<\varepsilon_0$ (see~\cite{AM}).
Then, using the equivalence between this metric and the Hausdorff distance on the set $\mathcal K$, see Appendix~\ref{appendix}, the following relation on the data $\mathcal D_{\mathcal H}(\partial \Omega_0,\partial \Omega_\varepsilon)=\mathcal O(\eps)$ now implies~\eqref{thm0}.\\

To obtain better estimate,  we introduce the following problem: Given a constant~$B_0$, find a domain~$A_\eps^0$ containing~$\Omega_0$ and function~$u_\eps^0$ defined on~$A_\eps^0 \setminus \overline{\Omega_0}$ such that
\begin{equation}\label{pbber1}
\left\{
\begin{array}{rcll}
\Delta u_\eps^0&=&0&\mbox{ in }A_\eps^0\setminus \overline{\Omega_0}, \\
u_\eps^0&=&1&\mbox{ on }\partial\Omega_0, \\
u_\eps^0&=&-\eps B_0&\mbox{ on }\partial A_\eps^0,\\
|\nabla u_\eps^0|&=&\lambda&\mbox{ on }\partial A_\eps^0.
\end{array}
\right.
\end{equation}
The main result presented in this paper is the following:
\begin{theorem}
\label{mainthm}
There exists $B_0\in \R$ such that the solution~$(A_\eps^0,u_\eps^0)$ of the Bernoulli problem~\eqref{pbber1} approaches the solution~$(A_\eps,u_\eps)$ of the Bernoulli problem~\eqref{pbberpert} at order two:
\begin{equation}
\label{thm}
\mathcal D_{\mathcal H}(A_\eps^0,A_\eps)=\mathcal O(\eps^2).
\end{equation}
\end{theorem}
For explicit calculations, the constant $B_0$ can be determined as follows.
We introduce the cell domain
\begin{equation}
\omega = \big\{ (R,\Theta) \in \R \times \T ~;~ R> - h(\Theta) \big\}.
\end{equation}
We solve the Laplace problem: Find a function $\widetilde u_1$ defined on~$\omega$, periodic with respect to~$\Theta$, with $\widetilde \nabla \widetilde u_1 \in L^2(\omega)$ and such that
\begin{equation}\label{pbcell1}
\left\{
\begin{array}{rcll}
\widetilde \Delta \widetilde u_1 &=& 0 &\mbox{ in } \omega, \\
\widetilde u_1(-h(\Theta),\Theta) &=& - h(\Theta) &\mbox{ on } \partial \omega.
\end{array}
\right.
\end{equation}
For such a function defined on the domain~$\omega$ the Laplacian operator reads:
$$
\widetilde \Delta = \frac{\partial^2}{\partial R^2} + \frac{\partial^2}{\partial \Theta^2}.
$$
The constant $B_0$ is then given by
\begin{equation}
\label{B0def}
B_0 =  \lambda \, \rho_0 \, \int_\T \widetilde u_1(0,\Theta) \, \mathrm d\Theta,
\end{equation}
the constant~$\rho_0$ corresponds to the radius of the domain~$A_0$ solution of the simple Bernoulli problem~\eqref{pbber0}.
\begin{remark}\label{rem2}
One of the key point of the proof is an estimate of kind:
\begin{equation}
\| u_\eps - ( u_\eps^0 + \eps \, \widetilde u_\eps^1) \|_{L^{\infty}(A_\eps\setminus \overline{\Omega_\eps})}= \mathcal O(\eps^2),
\end{equation}
for an "oscillating" function~$\widetilde u_\eps^1$ which exponentially decreases to the constant~$B_0$ far from the boundary~$\partial \Omega_\eps$.
\end{remark}

\section{Proofs}\label{proof}

\subsection{Well posedness}

First we prove that all the problems introduced in section \ref{stat} have a unique solution.

\subsubsection{Bernoulli problems}

The problems~\eqref{pbber0} and~\eqref{pbber1} are Bernoulli problems where the given domain is the unit ball $\Omega_0=\mathcal B(0,1)$. Consequently, their solution can be explicitly given (see for instance~\cite{FR}):\\
$\bullet$ The solution~$(A_0,u_0)$ of the problem~\eqref{pbber0} is invariant under rotation. It writes $A_0=\mathcal B(0,\rho_0)$ where $\rho_0>1$ satisfies:
\begin{equation}\label{expression-rho0}
\lambda \, \rho_0\ln \rho_0=1,
\end{equation}
and the function~$u_0$ is defined using polar coordinates by:
\begin{equation}\label{expression-u0}
\displaystyle u_0(r,\theta)=\frac{\ln(\rho_0)-\ln(r)}{\ln(\rho_0)}.
\end{equation}
$\bullet$ In the same way, we have $A_\eps^0=\mathcal B(0,\rho_\eps^0)$ where the constant~$\rho_\eps^0$ satisfies
\begin{equation}\label{expression-rhoeps}
\lambda \, \rho_\eps^0\ln \rho_\eps^0=1+B_0\eps.
\end{equation}
The explicit expression for~$u_\eps^0$ is then given by:
\begin{equation}\label{expression-ueps0}
\displaystyle u_\eps^0(r,\theta)=(1+B_0\eps)\frac{\ln(\rho_\eps^0)-\ln(r)}{\ln(\rho_\eps^0)}-B_0\eps.
\end{equation}
$\bullet$ The domain~$\Omega_\eps$ is defined using a graph of a lipschitz function (more precisely the function $\theta \mapsto 1-\eps h(\theta/\eps)$). Consequently this domain is starlike with respect to all points in an open ball centered at the origin.
We have the following existence result, proved in~\cite{AM}:
\begin{proposition}
The Bernoulli problem~\eqref{pbberpert} has a unique solution.
Moreover, the given domain~$\Omega_\eps$ being starlike with respect to all points in an open ball, the domain solution~$A_\eps$ is also starlike with respect to all points the same open ball.
\end{proposition}
As a consequence of this result, the boundary of the domain~$A_\eps$ can be writen as the graph of a lipschitz function:
\begin{equation}
\partial A_\eps =\{(r,\theta)\in \R_+ \times \T ~;~ r=\rho_\eps(\theta)\}.
\end{equation}

\subsubsection{Cell problem}

For the cell problem~\eqref{pbcell1}, an existence result of~$\widetilde u_1$ is given for instance in~\cite{NNRM}. 
Moreover we have the following  result:
\begin{proposition}\label{propcell1}
For all $(R,\Theta)\in \R_+ \times \T$ we have
\begin{equation}
\widetilde u_1(R,\Theta) = \int_\T \widetilde u_1(0,\Theta)\, \mathrm d\Theta + \sum_{k\in \mathbb Z \setminus \{0\}} \alpha_k \, \mathrm e^{-|k| R} \, \mathrm e^{ik\Theta},
\end{equation}
where each $\alpha_k$ is a constant.
\end{proposition}
{\bf Proof - }
We use the Fourier decomposition of the solution~$\widetilde u_1$ with respect to the periodic variable~$\Theta$. The Laplace equation results in an ordinary differential equation on each Fourier coefficient.
The proposition is then a consequence of the condition $\widetilde \nabla \widetilde u_1 \in L^2(\omega)$.
\begin{corollary}\label{corcell1}
For all $\mu < 1$ there exists $C(\mu) \geq 0$ such that for all $(R,\Theta)\in \R_+ \times \T$ we have:
\begin{equation}
\begin{aligned}
& \Big| \widetilde u_1(R,\Theta) - \int_\T \widetilde u_1(0,\Theta)\, \mathrm d\Theta \Big|õ'‰…あ\leq C(\mu) \, \mathrm e^{-\mu R},\\
& \Big| \widetilde \nabla \widetilde u_1(R,\Theta) \Big|õ'‰…あ\leq C(\mu) \, \mathrm e^{-\mu R}.
\end{aligned}
\end{equation}
\end{corollary}
%

\subsection{Laplace system with oblique boundary conditions for the error}

Now we prove the estimate announced in Remark~\ref{rem2}.

\subsubsection{Building of the function $\widetilde u_\eps^1$}~\par

$\bullet$ Besides the introduction of the function~$\widetilde u_1$, we introduce by induction the following functions as solution of Laplace problems on the cell domain~$\omega$:
The functions~$\widetilde u_j$, $j\geq 2$, periodic with respect to~$\Theta$, with $\widetilde \nabla \widetilde u_j \in L^2(\omega)$ and such that:
\begin{equation}\label{pbcellj}
\left\{
\begin{aligned}
& \widetilde \Delta \widetilde u_j = \sum_{k=1}^{j-1} (-1)^{j-k} R^{j-k-1} \left( \frac{\partial \widetilde u_{k}}{\partial R} - (j-k+1) R \, \frac{\partial^2 \widetilde u_{k}}{\partial \Theta^2} \right) & \mbox{in } \omega, \\
& \widetilde u_j(-h(\Theta),\Theta) = 0 & \mbox{on } \partial \omega.
\end{aligned}
\right.
\end{equation}
Note that the behavior of these functions~$\widetilde u_j$, $j\geq 2$, for large values of~$R$ is similar to those of the function~$\widetilde u_1$ (see Proposition~\ref{propcell1} and Corollary~\ref{corcell1}).
By induction on~$j$ we prove that the right hand side members of the Laplace equation~\eqref{pbcellj} are free average with respect to the periodic variable~$\Theta$ and, like in the proof of the proposition~\ref{propcell1}, see also the same kind of proof in~\cite{ChupinMartin12}, the following results:
\begin{proposition}\label{propcellj}
For each $j\geq 2$ there exists a solution~$\widetilde u_j$ of the Laplace problem~\eqref{pbcellj}.
Moreover for $(R,\Theta)\in \R_+ \times \T$ we have
\begin{equation}
\widetilde u_j(R,\Theta) = \int_\T \widetilde u_j(0,\Theta)\, \mathrm d\Theta + \sum_{k\in \mathbb Z \setminus \{0\}} \alpha_{k,j}(R) \, \mathrm e^{-|k| R} \, \mathrm e^{ik\Theta},
\end{equation}
where each $\alpha_{k,j}$ is a polynomial function.
In particular for all $\mu < 1$ there exists $C(\mu) \geq 0$ such that for $(R,\Theta)\in \R_+ \times \T$ we have
\begin{equation}
\begin{aligned}
& \Big| \widetilde u_j(R,\Theta) - \int_\T \widetilde u_j(0,\Theta)\, \mathrm d\Theta \big| \leq C(\mu) \, \mathrm e^{-\mu R},\\
& \Big| \widetilde \nabla \widetilde u_j(R,\Theta) \Big| \leq C(\mu) \, \mathrm e^{-\mu R}.
\end{aligned}
\end{equation}
\end{proposition}
$\bullet$ We introduce
\begin{equation}\label{construction-u1tilde}
\widetilde u_\eps^1(r,\theta) = \lambda \, \rho_0 \sum_{j=1}^J \eps^{j-1} \widetilde u_j \Big(\frac{r-1}{\eps},\frac{\theta}{\eps} \Big),
\end{equation}
the value of the integer~$J$ will be later selected.
Due to the change of variable $\displaystyle R=\frac{r-1}{\eps}$, $\displaystyle \Theta=\frac{\theta}{\eps}$ we have for each $j\geq 1$:
\begin{align*}
\Delta \left[  \widetilde u_j \Big(\frac{r-1}{\eps},\frac{\theta}{\eps} \Big) \right]
= \frac{1}{\eps^2} \frac{\partial^2 \widetilde u_j}{\partial R^2} (R,\Theta)
&   + \frac{1}{\eps (1+\eps R)} \frac{\partial \widetilde u_j}{\partial R} (R,\Theta)\\
&   + \frac{1}{\eps^2 (1+\eps R)^2} \frac{\partial^2 \widetilde u_j}{\partial \Theta^2} (R,\Theta).
\end{align*}
We deduce the following asymptotic development with respect to~$\eps$:
\begin{align*}
\Delta \left[  \widetilde u_j \Big(\frac{r-1}{\eps},\frac{\theta}{\eps} \Big) \right]
& = \frac{1}{\eps^2} \widetilde \Delta \widetilde u_j (R,\Theta) \\
+ \sum_{k=0}^{\infty} & \eps^{k-1} (-1)^k R^k \left( \frac{\partial \widetilde u_j}{\partial R}(R,\Theta) - (k+2) R \frac{\partial^2 \widetilde u_j}{\partial \Theta^2}(R,\Theta) \right).
\end{align*}
$\bullet$
By construction of~$\widetilde u_\eps^1$, see its expression~\eqref{construction-u1tilde}, we deduce that $\Delta \widetilde u_\eps^1 = \eps^{J-2} F_\eps$, where $F_\eps(r,\theta) = \mathcal O(1)$ with respect to~$\eps$.
Note that in the sequel we will choose $J=3$ so that
\begin{equation}\label{Delta-utilde}
\Delta \widetilde u_\eps^1 = \mathcal O(\eps).
\end{equation}

$\bullet$
Note also that the value of the function~$\widetilde u_\eps^1$ on the boundary~$\partial \Omega_\eps$ is given by
\begin{equation}\label{boundary1-utilde}
\widetilde u_\eps^1\vert_{\partial \Omega_\eps} = - \lambda \, \rho_0 h\left(\frac{\cdot}{\eps}\right).
\end{equation}
$\bullet$
On the boundary~$\partial A_\eps$ parametrised by the function~$\rho_\eps$ we have, for all $\theta\in \T$
\begin{align*}
\widetilde u_\eps^1 (\rho_\eps(\theta),\theta)  - \lambda \, \rho_0 & \int_\T \widetilde u_1(0,\Theta)\, \mathrm d\Theta \\
= & \lambda \, \rho_0 \left( \widetilde u_1\Big(\frac{\rho_\eps(\theta)-1}{\eps},\frac{\theta}{\eps}\Big) - \int_\T \widetilde u_1(0,\Theta)\, \mathrm d\Theta \right) \\
& + \eps \lambda \, \rho_0 \sum_{j=2}^J \eps^{j-2} \widetilde u_j\Big(\frac{\rho_\eps(\theta)-1}{\eps},\frac{\theta}{\eps}\Big).
\end{align*}
Using Corollary~\ref{corcell1}, Proposition~\ref{propcellj} and the fact that for all $\theta\in \T$ we have $\rho_\eps(\theta)-1=\mathcal O(1)$ (since $\displaystyle \lim_{\eps\to 0} \rho_\eps(\theta) = \rho_0$ for all $\theta \in \T$, see~\cite{AM}) we deduce
\begin{equation}\label{boundary2-utilde}
\Big| \widetilde u_\eps^1\vert_{\partial A_\eps}  - \lambda \, \rho_0 \int_\T \widetilde u_1(0,\Theta)\, \mathrm d\Theta \Big| = \mathcal O(\eps).
\end{equation}
$\bullet$
In the same way, using the estimate for the gradients in Corollary~\ref{corcell1} and Proposition~\ref{propcellj}, we have for instance
\begin{equation}\label{boundary3-utilde}
\nabla \widetilde u_\eps^1\vert_{\partial A_\eps} = o(\eps^k)\qquad \text{for all $k\in \mathbb N$}.
\end{equation}

\subsubsection{System satisfied by the error}

We introduce the difference $v=u_\eps - (u_\eps^0+ \eps \widetilde u_\eps^1)$.
We will prove that $v$ satisfies a Laplace problem on the domain~$A_\eps \setminus \overline{\Omega_\eps}$ with oblique boundary conditions on the exterior boundary~$\partial A_\eps$.\\

$\bullet$ We are first interested in the value of~$\Delta v$ on $A_\eps\setminus \overline{\Omega_\eps}$.\\
By definition of~$u_\eps$ (which satisfies the Bernoulli problem~\eqref{pbberpert}), we have $\Delta u_\eps = 0$ on $A_\eps\setminus \overline{\Omega_\eps}$.
The function $u_\eps^0$, solution of the Bernoulli problem~\eqref{pbber1}, is {\it a priori} only defined on~$A_\eps^0\setminus \overline{\Omega_0}$. But its explicit expression given by~\eqref{expression-ueps0} can be extended on~$A_\eps\setminus \overline{\Omega_\eps}$, by preserving the relation~$\Delta u_\eps^0 = 0$.
Finally, the approximation~\eqref{Delta-utilde} on the oscillating contribution~$\Delta \widetilde u_\eps^1$ implies that
\begin{equation}\label{Delta-v}
\Delta v = \mathcal O(\eps^2).
\end{equation}

$\bullet$ The second step consists in obtaining the value of~$v$ on~$\partial \Omega_\eps$.\\
From the Bernoulli problem~\eqref{pbberpert}, we know that $u_\eps\vert_{\partial \Omega_\eps} = 1$.\\
The Taylor development of the function~$u_\eps^0$ (or more precisely its extension) implies that for all~$\theta\in \T$ we have
$$u_\eps^0(1-\eps h\left(\frac{\theta}{\eps}\right),\theta) = u_\eps^0(1,\theta) - \eps h\left(\frac{\theta}{\eps}\right) \partial_r u_\eps^0(1,\theta) + \mathcal O(\eps^2).$$
From the analytic expression of~$u_\eps^0$, see~\eqref{expression-ueps0}, we deduce that for all~$\theta\in \T$ we have
 $$u_\eps^0(1-\eps h\left(\frac{\theta}{\eps}\right),\theta) = 1 + \eps \lambda \rho_0 h\left(\frac{\theta}{\eps}\right) + \mathcal O(\eps^2).$$
Note that we have used the relation $\rho_\eps^0 = \rho_0 + \mathcal O(\eps)$ which directly follows from the expressions of $\rho_\eps^0$ and~$\rho_0$ given by~\eqref{expression-ueps0} and~\eqref{expression-u0} respectively.
Using the estimate~\eqref{boundary1-utilde} of~$\widetilde u_\eps^1$ on the boundary~$\partial \Omega_\eps$ we deduce that
\begin{equation}\label{boundary1-v}
v\vert_{\partial \Omega_\eps} =  \mathcal O(\eps^2).
\end{equation}

$\bullet$ We now compute the value of~$v$ on the boundary~$\partial A_\eps$.\\
By definition of the solution~$u_\eps$ of the problem~\eqref{pbberpert}, we have
$u_\eps \vert_{\partial A_\eps} = 0$.\\
Recall that the boundary $\partial A_\eps$ is parametrised using the function~$\rho_\eps$ whereas the boundary~$\partial A_\eps^0$ is the ball $\mathcal B(0,\rho_\eps^0)$.
We deduce that for all $\theta\in \T$ we have
$$
u_\eps^0(\rho_\eps(\theta),\theta) = u_\eps^0(\rho_\eps^0,\theta) + ( \rho_\eps(\theta) - \rho_\eps^0 ) \partial_r u_\eps^0(\rho_\eps^0,\theta) + \mathcal O((\rho_\eps(\theta) - \rho_\eps^0)^2).
$$
Since the function~$u_\eps^0$ solves the Bernoulli problem~\eqref{pbber1}, using the fact that~$\partial_r$ exactly corresponds to the normal derivative on the circle~$\partial A_\eps^0$, and using that for all $\theta\in \T$ we have $\rho_\eps(\theta) - \rho_\eps^0 = \mathcal O( \eps)$, we obtain
$$
u_\eps^0(\rho_\eps(\theta),\theta) = -\eps\, B_0 - ( \rho_\eps(\theta) - \rho_\eps^0 ) \lambda + \mathcal O(\eps^2).
$$
Choosing $B_0=\lambda \rho_0 \int_\T \widetilde u_1(0,\Theta) \, \mathrm d\Theta$ , we can then use the result~\eqref{boundary2-utilde} on the oscillating term~$\widetilde u_\eps^1$ to write
\begin{equation}\label{dirirho}
v\vert_{\partial A_\eps} = ( \rho_\eps - \rho_\eps^0 ) \lambda + \mathcal O(\eps^2).
\end{equation}

$\bullet$ In the same way, we obtain the value of some derivative of~$v$ on~$\partial A_\eps$:\\
By definition of~$u_\eps$, and denoting by~$n_\eps$ the outward unitary normal to the boundary~$\partial A_\eps$, we have 
\begin{equation}
\label{clv1}
\nabla u_\eps \vert_{\partial A_\eps} = -\lambda n_\eps.
\end{equation}
Using the Taylor formulae, we obtain
\begin{equation}
\label{clv2-0}
\nabla u_\eps^0 \vert_{\partial A_\eps} = \nabla u_\eps^0 \vert_{\partial A_\eps^0} + (\rho_\eps - \rho_\eps^0) (\partial_r \nabla u_\eps^0) \vert_{\partial A_\eps^0} + \mathcal O((\rho_\eps - \rho_\eps^0)^2).
\end{equation}
Due to the Bernoulli problem~\eqref{pbber1} satisfied by~$u_\eps^0$, we have $$\nabla u_\eps^0 \vert_{\partial A_\eps^0} = -\lambda n_\eps^0,$$ where~$n_\eps^0$ denotes the outward unitary normal to the boundary~$\partial A_\eps^0$.\\
Moreover the expression of the laplacian operator in polar coordinates on the boundary~$\partial A_\eps^0$ reads: 
$$\Delta u_\eps^0 = \frac{1}{(\rho_\eps^0)^2}\partial_\theta^2 u_\eps^0 + \frac{1}{\rho_\eps^0} \partial_{n_\eps^0} u_\eps^0 + \partial^2_{n_\eps^0} u_\eps^0.$$ 
Since~$u_\eps^0$ solves the Bernoulli problem~\eqref{pbber1}, we deduce that $\displaystyle(\partial_r \nabla u_\eps^0) \vert_{\partial A_\eps^0} = \frac{\lambda}{\rho_\eps^0}n_\eps^0$.\\
The equality~\eqref{clv2-0} now reads
\begin{equation}
\label{clv2}
\nabla u_\eps^0 \vert_{\partial A_\eps} = - \lambda n_\eps^0 + \frac{\rho_\eps - \rho_\eps^0}{\rho_\eps^0}\lambda n_\eps^0 + \mathcal O(\eps^2).
\end{equation}
Substracting \eqref{clv2} from \eqref{clv1}, and noting that the oscillating contribution is not dominating (see the equation~\eqref{boundary3-utilde}) we get
$$
\nabla v\vert_{\partial A_\eps} = - \lambda (n_\eps - n_\eps^0) - \frac{\rho_\eps - \rho_\eps^0}{\rho_\eps^0}\lambda n_\eps^0 + \mathcal O(\eps^2).
$$
Taking the scalar product with $\widetilde n=\left(n_\eps + n_\eps^0\right)/\|n_\eps + n_\eps^0\|$ gives:
\begin{equation}
\label{neumrho}
\partial_{\widetilde n} v\vert_{\partial A_\eps} = - \frac{\rho_\eps - \rho_\eps^0}{\rho_\eps^0}\lambda \left( n_\eps^0 \cdot \widetilde n \right) + \mathcal O(\eps^2).
\end{equation}
$\bullet$ We finally deduce the system satisfy by~$v$:
\begin{equation}\label{pbv}
\left\{
\begin{array}{rcll}
\Delta v&=&f& \mbox{ in } A_\eps\setminus \overline{\Omega_\eps}, \\
v&=&g& \mbox{ on } \partial\Omega_\eps, \\
\partial_{\widetilde n} v + \frac{1}{\rho_\eps^0} \left(n_\eps^0 \cdot \widetilde n\right) \, v&=&h& \mbox{ on } \partial A_\eps,
\end{array}
\right.
\end{equation}
with $f=\mathcal O(\eps^2)$, $g=\mathcal O(\eps^2)$ and $h=\mathcal O(\eps^2)$. The last condition is an oblique boundary condition obtained combining~\eqref{dirirho} and~\eqref{neumrho}.


\newpage
\subsection{Estimate of the error}
\subsubsection{Bounds on a Laplace type problem}

In this section, we show that the regular function $v$ satisfying~\eqref{pbv} on the regular domain~$A_\eps\setminus \overline{\Omega_\eps}$ satisfies $v=\mathcal O(\eps^2)$.\\

Let $\gamma=\frac{n_\eps^0 \cdot \widetilde n}{\rho_\eps^0}$. It is clear that $\gamma>\sqrt2/(2\rho_\varepsilon^0)>\sqrt2/(2\rho_0)>0$, thus one can choose a positive number $R$ which does not depend on $\varepsilon$ such that the function $\varphi(x,y) = R^2-\frac{1}{4}(x^2+y^2)$ satisfies:
\begin{equation}\label{pbv-varphi}
\left\{
\begin{array}{rcll}
- \Delta \varphi&=& 1& \mbox{ in } A_\eps\setminus \overline{\Omega_\eps}, \\
\varphi&\geq&1& \mbox{ on } \partial\Omega_\eps,\\
\displaystyle \alpha\partial_{n_\eps}\varphi+\beta\partial_{\tau_\eps}\varphi+\gamma\varphi&\geq&1& \mbox{ on } \partial A_\eps,
\end{array}
\right.
\end{equation}
where the coefficients $\alpha=\tilde n\cdot{n_\eps}$ and $\beta=\tilde n-\alpha{n_\eps} \in \R$  are chosen so that the boundary condition on $\partial A_\eps$ corresponds to the boundary condition of system \eqref{pbv}, $\partial_{\tau_\eps}\varphi\in\R$ denoting the tangential derivative on $\partial A_\eps$.

Note that $\sqrt2/2<\alpha\le 1$.

Consequently the function $V=\max\{|f|_\infty, |g|_\infty, |h|_\infty\} \varphi - v$ is a regular function which satisfies
\begin{equation}\label{pbv-varphi-v}
\left\{
\begin{array}{rcll}
- \Delta V&\geq& 0& \mbox{ in } A_\eps\setminus \overline{\Omega_\eps}, \\
V&\geq&0& \mbox{ on } \partial\Omega_\eps, \\
\displaystyle \alpha\partial_{n_\eps}V+\beta\cdot\partial_{\tau_\eps}V+\gamma V&\geq&0& \mbox{ on } \partial A_\eps.
\end{array}
\right.
\end{equation}

Moreover, the domain~$A_\eps\setminus \overline{\Omega_\eps}$ is regular since, by hypothesis, $\Omega_\eps$ is regular and since the domain~$A_\eps$ is regular too (see \cite{AM}).
This allows us to deduce the following result:

\begin{lemma}
\label{monoto}
The function $V$ is nonnegative on $A_\eps\setminus \overline{\Omega_\eps}$.
\end{lemma}
{\bf Proof of Lemma~\ref{monoto} - }
Since $\Delta V\le0$ the minimum of~$V$ in $\overline{A_\eps\setminus \overline{\Omega_\eps}}$ is achieved on  its boundary $\partial{(A_\eps\setminus \overline{\Omega_\eps})}=\partial{A_\eps}\cup\partial{\Omega_\eps}$.
If it is achieved on $\partial\Omega_\eps$, then the proof is completed.
If it is achieved on $y\in\partial{A_\eps}$, then, since $V$ is regular,  we have $\partial_{n_\eps}V(y)\le0$ and $\partial_{\tau_\eps}V(y)=0$. Since $\alpha\geq 0$ we deduce that $V(y)=\mathcal \gamma^{-1}\left(h(y)- \alpha(y)\partial_{n_\eps}V(y)-\beta(y)\cdot\partial_{\tau_\eps}V(y)\right)\ge0$.\\

Lemma~\ref{monoto} implies that $v\leq\max\{|f|_\infty, |g|_\infty, |h|_\infty\} |\varphi|_\infty$. 

In the same way, using $\tilde V=-\max\{|f|_\infty, |g|_\infty, |h|_\infty\} \varphi - v$ instead of~$V$ we show that $\tilde V\leq 0$.
We conclude that
\begin{equation} 
\label{Maj}
|v|_\infty \leq \max\{|f|_\infty, |g|_\infty, |h|_\infty\} |\varphi|_\infty=\mathcal O(\varepsilon^2).
\end{equation}

Due to the equation~\eqref{dirirho} we deduce that $\rho_\eps - \rho_\eps^0 = \mathcal O(\eps^2)$, which is equivalent to $\mathcal D_{\mathcal H}(A_\eps^0,A_\eps)=\mathcal O(\eps^2)$.
Theorem~\ref{mainthm} is proved.

\section{Numerical results}\label{numerics}

This section is devoted to numerical results to check the conclusion of Theorem \ref{mainthm}, first by comparing the solution of \eqref{pbberpert} and \eqref{pbber1} and then by observing that the roughness on the free boundary are much smaller than those of the fixed boundary.\\

Two kinds of perturbations have been considered, which are the $2\pi$-periodic functions denoted by~$h_1$ and~$h_2$ defined on $(0,2\pi)$ by:
\begin{itemize}
\item $h_1(\alpha)= 1-\cos(\alpha)$,
\item $h_2(\alpha)= \pi-|\alpha-\pi|$.
\end{itemize}

\subsection{Numerical approximation for the oscillating problem~\eqref{pbber1}}

In this section, we effectively compute the value of the constant~$B_0$ which appears in the Bernoulli problem~\eqref{pbber1}.
In practice, we need to solve the problem~\eqref{pbcell1} which is defined in an infinite domain~$\omega$, and then evaluate the integral~$\int_\T \widetilde u_1(0,\Theta) \, \mathrm d\Theta$ to deduce the value of the constant~$B_0$ {\it via} the formulae~\eqref{B0def}. For a numerical approach we introduce the ``truncated'' cell domain
\begin{equation}
\omega_M = \big\{ (R,\Theta) \in \R_+ \times \T ~;~ M>R> - h(\Theta) \big\},
\end{equation}
where~$M$ has been chosen large, and solve the Laplace problem: Find a function $\widetilde u_1^M$ defined on~$\omega_M$, periodic with respect to~$\Theta$ and such that
\begin{equation}\label{pbcell2}
\left\{
\begin{array}{rcll}
\widetilde \Delta \widetilde u_1^M &=& 0 &\mbox{ in } \omega_M, \\
\widetilde u_1^M(-h(\Theta),\Theta) &=& - h(\Theta) &\mbox{ for } \Theta\in\T, \\
\partial_R\widetilde u_1^M(M,\Theta)&=&0&\mbox{ for } \Theta\in\T.
\end{array}
\right.
\end{equation}
The solution of this problem being exponentially close to the solution of~\eqref{pbcell1} for large $M$, we have:
\begin{equation}
\int_\T \widetilde u_1(0,\Theta) \, \mathrm d\Theta \approx \int_\T \widetilde u_1^M(0,\Theta) \, \mathrm d\Theta.
\end{equation}
The computations of these integrals were performed with the FreeFem++ program (This software, see http://www.freefem.org/ff++ is based on weak formulation of the problem and finite elements method). We choose $M=6$ and the mesh is composed of about $2\cdot 10^6$ triangles.
\par
These computations give:
\begin{itemize}
\item $\displaystyle\int_\T \widetilde u_1(0,\Theta) \, \mathrm d\Theta\approx-0.58738$
for the shape $h_1(\alpha)= 1-\cos(\alpha)$,

\item $\displaystyle\int_\T \widetilde u_1(0,\Theta) \, \mathrm d\Theta\approx-0.87754$
for the shape $h_2(\alpha)= \pi-|\alpha-\pi|$.
\end{itemize}

\subsection{Comparison of problems \eqref{pbberpert} and \eqref{pbber1}}\label{comp}

The aim of these first tests is to show that the solution of \eqref{pbber1} approaches the solution of \eqref{pbberpert} at second order with respect to~$\eps$.
Note that all these tests use the algorithm presented in~\cite{BCT} to solve the Bernoulli problem \eqref{pbberpert}.\\

Here, we take $\lambda=2\mathrm e^{-1/2}$ so that the domain solution without roughness is a disc of radius $\rho_0=\mathrm e^{1/2}$. Thus, we have $\lambda\rho_0=2$, and we then:
\begin{itemize}
\item $B_0=\lambda\rho_0\displaystyle\int_\T \widetilde u_1(0,\Theta) \, \mathrm d\Theta\approx-1.17476$ 
for the shape $h_1(\alpha)= 1-\cos(\alpha)$,

\item $B_0=\lambda\rho_0\displaystyle\int_\T \widetilde u_1(0,\Theta) \, \mathrm d\Theta\approx-1.75508$ 
for the shape $h_2(\alpha)= \pi-|\alpha-\pi|$.
\end{itemize}

We check that 
$$\frac{\mathcal D_{\mathcal H}(A_\varepsilon,A_\varepsilon^0)}{\varepsilon^2} ~ \text{is bounded}$$
with respect to~$\eps$ according to Theorem \ref{mainthm}: The computed values are reported in Table~\ref{tablecomp2} and~\ref{tablecomp3-b}.
On coarse meshes (see the two first lines of tables  \ref{tablecomp2} and \ref{tablecomp3-b}), some grid effects deteriorate the quality of numerical results for small $\varepsilon$ due to the small size of roughness.

The numerical simulation on the finer mesh (see last line of tables  \ref{tablecomp2} and \ref{tablecomp3-b}) seems accurate enough to confirm the conclusion of Theorem \ref{mainthm}.
\begin{table}[h!]
\begin{tabular}{|r|c|c|c|}
\hline
&$\varepsilon=0.1$&$\varepsilon=0.05$& $\varepsilon=0.025$\\
\hline
$\delta x=6\times10^{-3}$&0.1524&0.0186&0.2832 \\
$\delta x=3\times10^{-3}$&0.1551&0.0238&0.0451 \\
$\delta x=1.5\times10^{-3}$&0.1556&0.0290&0.0145 \\
\hline
\end{tabular}
\vspace{0.3cm}
\caption{Values of $\displaystyle \frac{\mathcal D_{\mathcal H}(A_\varepsilon,A_\varepsilon^0)}{\varepsilon^2}$ for the perturbation $h_1$.}
\label{tablecomp2}
\end{table}

\begin{table}[h!]
\begin{tabular}{|r|c|c|c|}
\hline
&$\varepsilon=0.1$&$\varepsilon=0.05$& $\varepsilon=0.025$\\
\hline
$\delta x=6\times10^{-3}$&0.2579&0.1382&1.265\phantom{3} \\
$\delta x=3\times10^{-3}$&0.2885&0.1399&0.2028 \\
$\delta x=1.5\times10^{-3}$&0.2850&0.1302&0.0362 \\
\hline
\end{tabular}
\vspace{0.3cm}
\caption{Values of $\displaystyle \frac{\mathcal D_{\mathcal H}(A_\varepsilon,A_\varepsilon^0)}{\varepsilon^2}$ for the perturbation $h_2$.}
\label{tablecomp3-b}
\end{table}

\subsection{The oscillations of the free boundary}

\begin{figure}[!ht]
\centerline{\includegraphics[width=7.5cm,height=7.5cm]{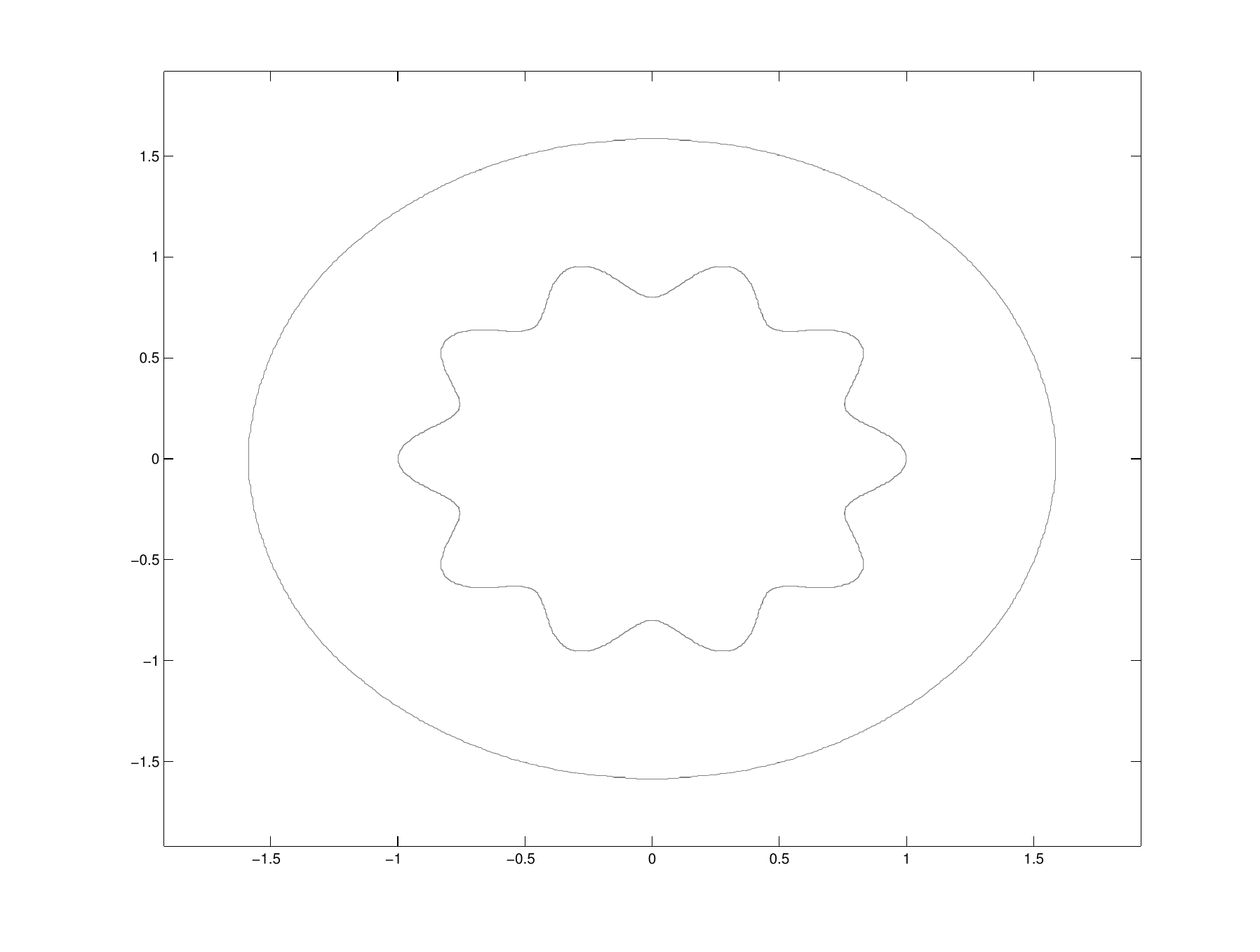}}
\caption{Given domain and its associated solution for the perturbation $h_1$
($\lambda=2\mathrm e^{-1/2}$, $\rho_0=\mathrm e^{1/2}$, $\varepsilon=0.1$).}
\label{fig_cas1}
\end{figure}

We now present some figures to show the oscillations observed on the computed solution. We plot on Figure \ref{fig_cas1} the given and solution domains corresponding to the first test case of the previous section. As indicated by table \ref{tablecomp2}, we just observe that the solution domain is very close to a disc.

In order to observe small oscillations on the solution domains, we changed the value of $\lambda$ so that $\lambda=8\mathrm e^{-1/8}$, and then $\rho_0=\mathrm e^{1/8}$ which means that the solution domain comes closer to the given domain which is (intuitively) a better situation to observe oscillations. 
The function $h$ corresponding to the given domain is the function $h_1$ of section \ref{comp}, and we take $\delta_x=6\times10^{-3}$.

\begin{figure}[!ht]
\centerline{\includegraphics[width=7.5cm,height=7.5cm]{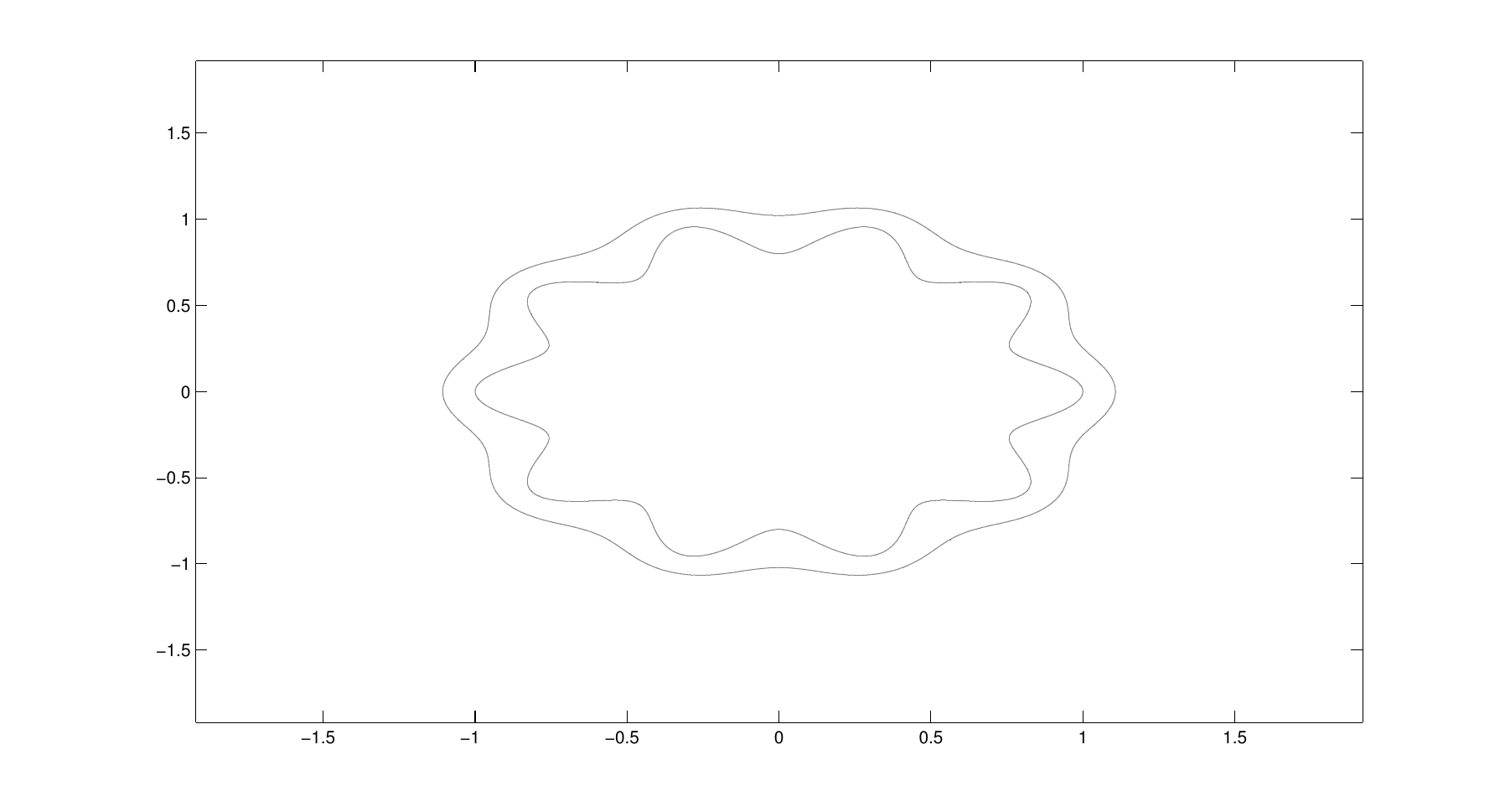}
\includegraphics[width=7.5cm,height=7.5cm]{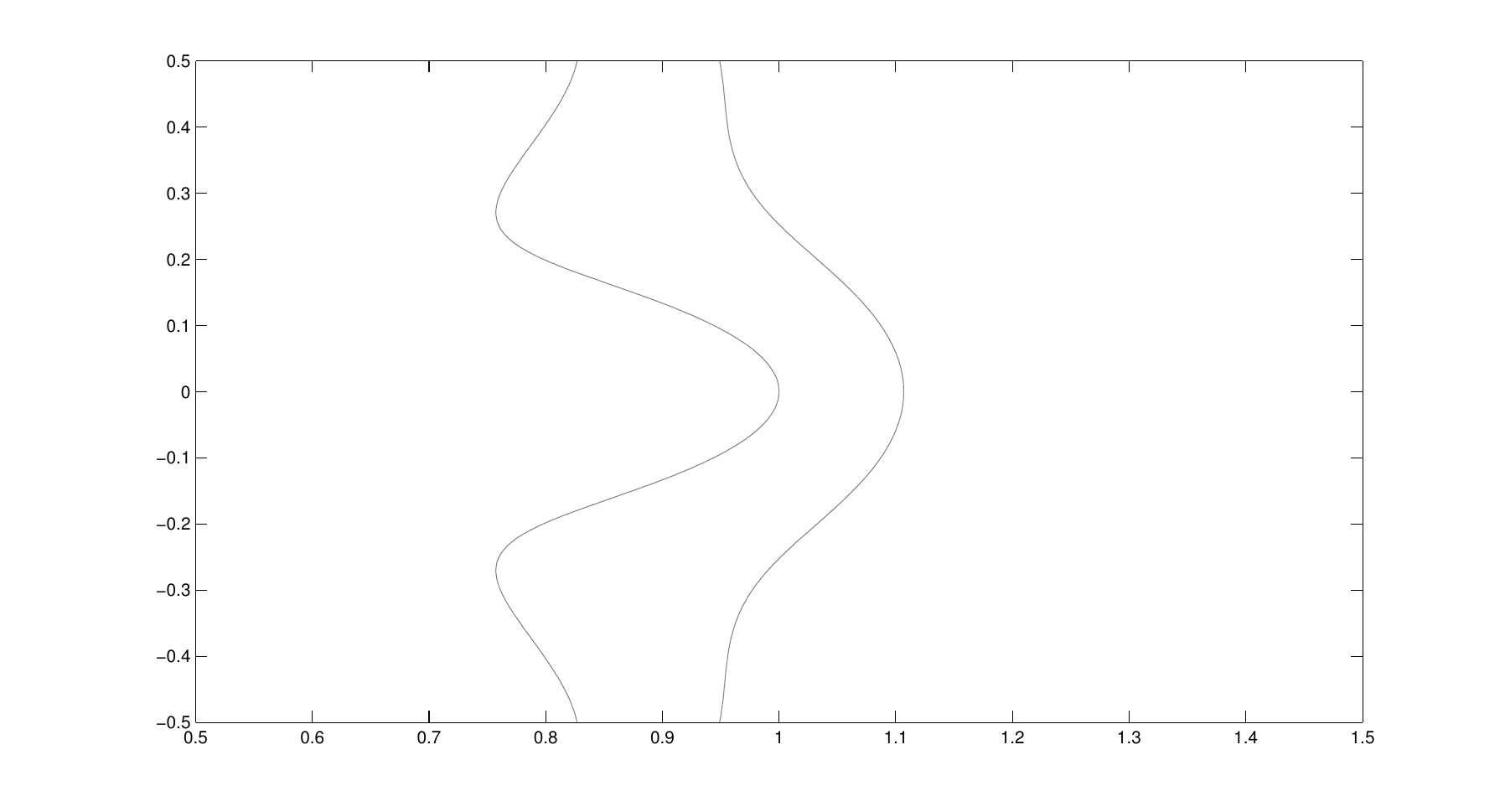}
}
\caption{Given domain and its associated solution for the perturbation $h_1$ with $\lambda=8\, \mathrm e^{-1/8}$, $\rho_0=\mathrm e^{1/8}$, $\varepsilon=0.1$.}
\label{fig_cas2}
\end{figure}

On Figure \ref{fig_cas2}, for $\varepsilon=0.1$, we can observe that the solution domain is not a circle and that some small oscillations seem to be induced by the oscillations of the given domain $\Omega_\varepsilon$.

\begin{figure}[!ht]
\centerline{\includegraphics[width=7.5cm,height=7.5cm]{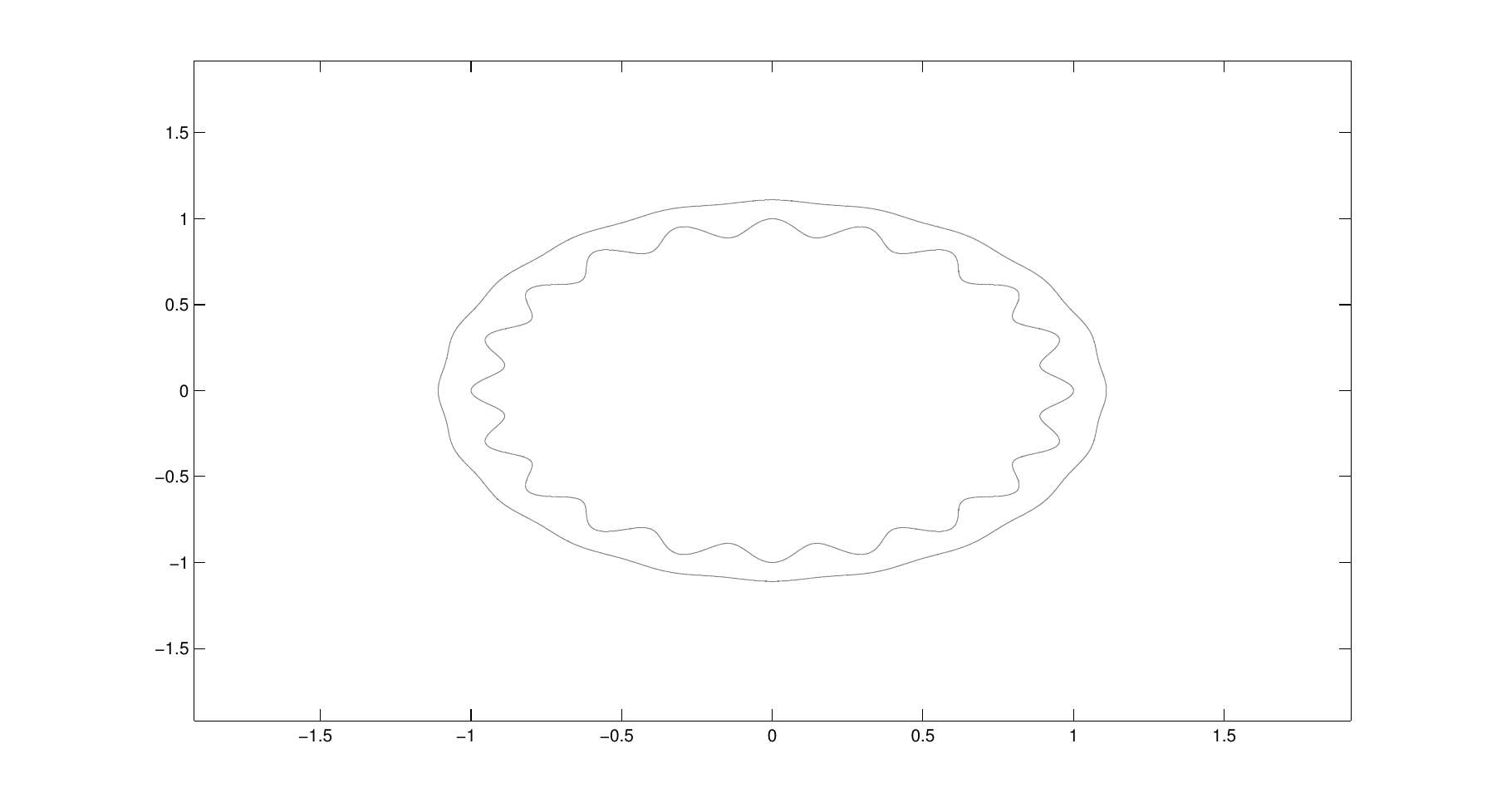}
\includegraphics[width=7.5cm,height=7.5cm]{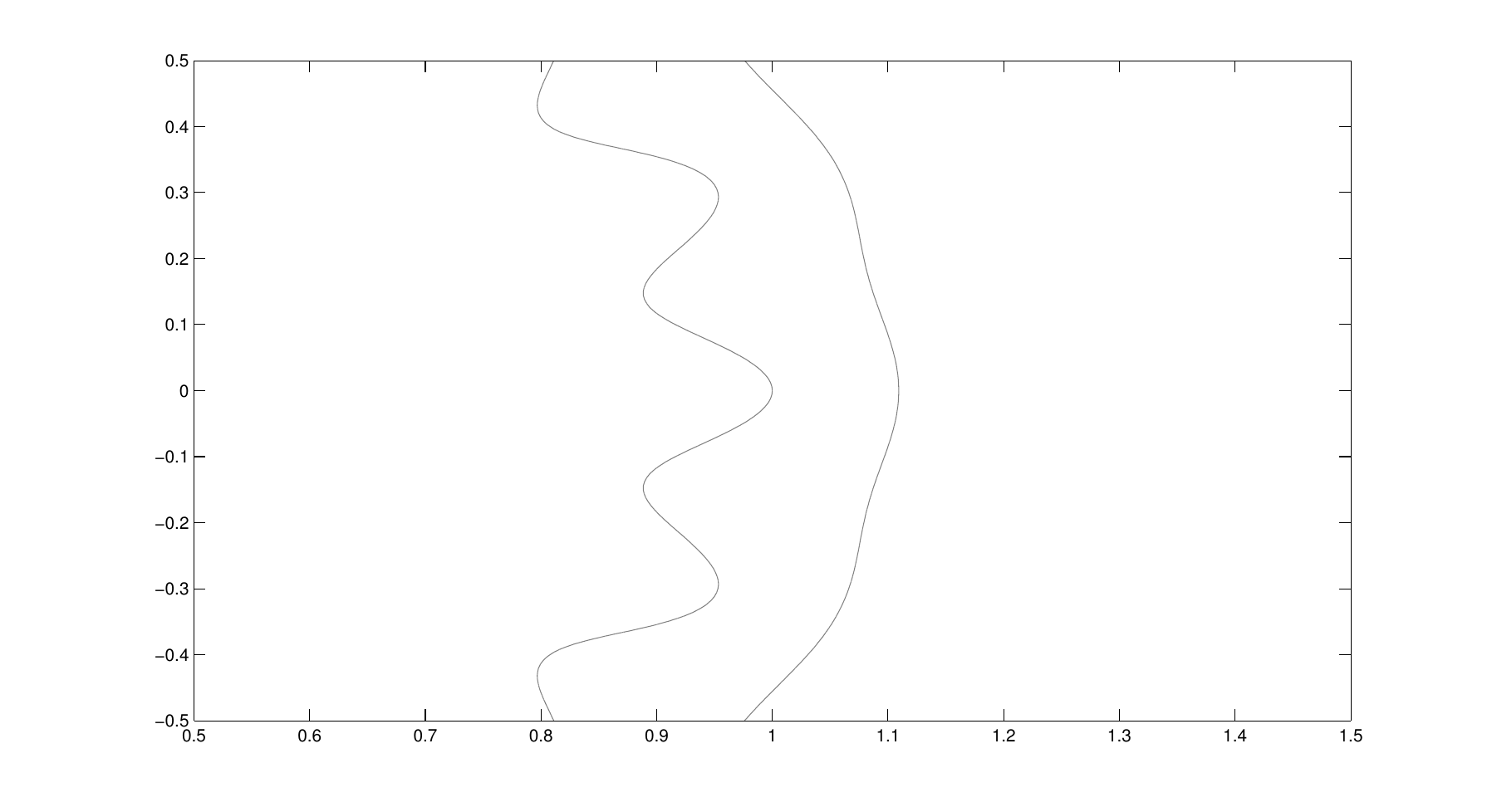}
}
\caption{Given domain and its associated solution for the perturbation $h_1$ with $\lambda=8\, \mathrm e^{-1/8}$, $\rho_0=\mathrm e^{1/8}$, $\varepsilon=0.05$.}
\label{fig_cas2a}
\end{figure}

On Figure \ref{fig_cas2a}, for a smaller value of $\varepsilon$, we need to zoom on a part of the picture to observe the same behaviour of the solution domain. 

\begin{figure}[!ht]
\centerline{\includegraphics[width=7.5cm,height=7.5cm]{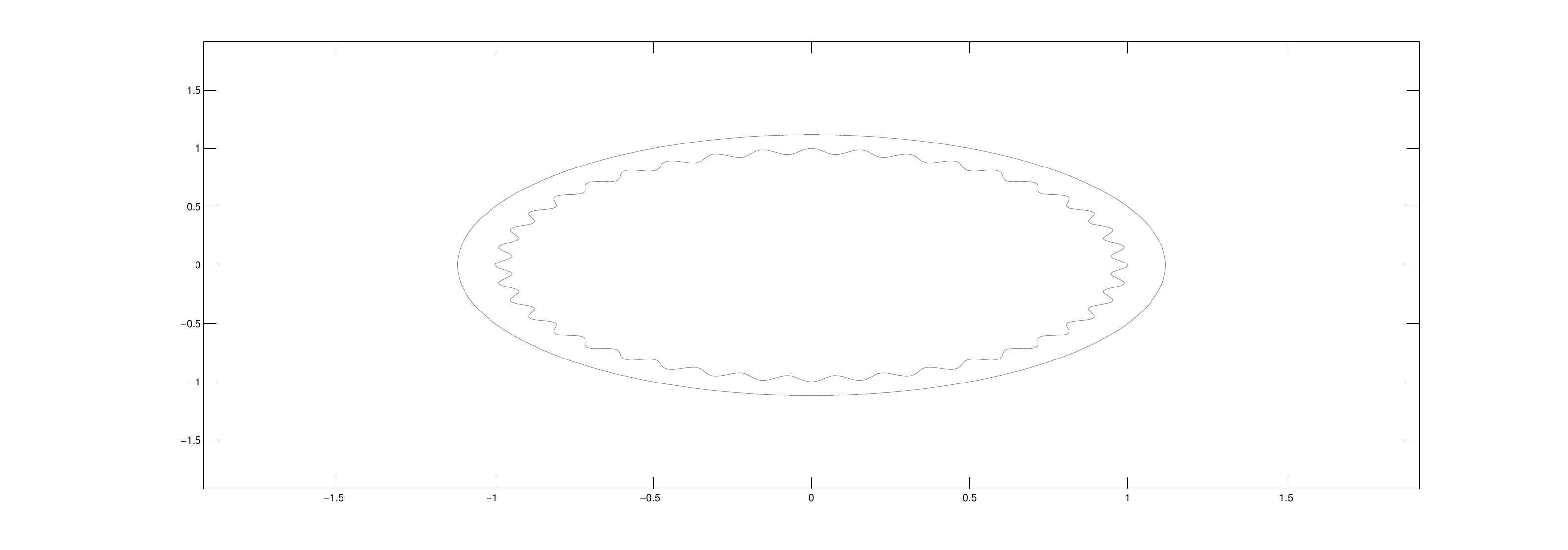}
\includegraphics[width=7.5cm,height=7.5cm]{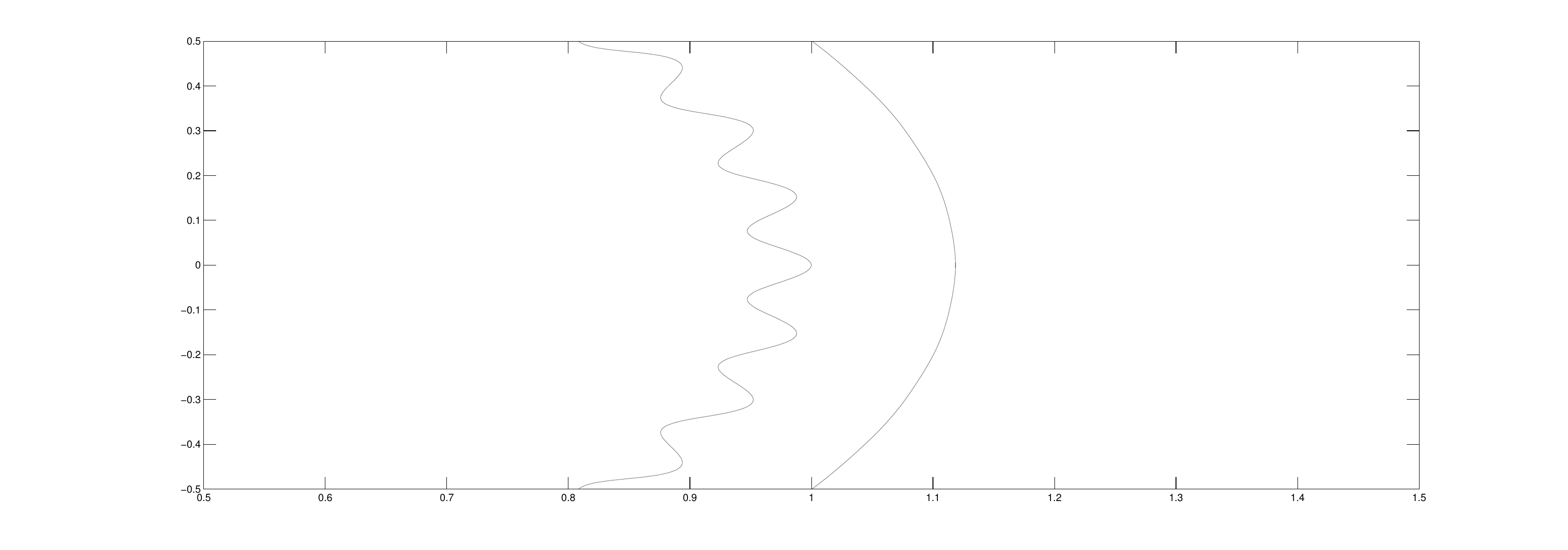}}
\caption{Given domain and its associated solution for the perturbation $h_1$ with $\lambda=8\mathrm e^{-1/8}$, $\rho_0=\mathrm e^{1/8}$, $\varepsilon=0.025$.}
\label{fig_cas2b}
\end{figure}

On Figure \ref{fig_cas2b}, the value of $\varepsilon$ is so small that the oscillation on the solution domain can no more be observed even on the zoomed figure. These tests also confirm the qualitative conclusion of Theorem \ref{mainthm}.


\section{Conclusion}\label{concl}

In this work, we have given a result concerning the solution of the free boundary Bernoulli problem with a given rough domain. In the case where the given domain is a disc with roughness, we show that the solution is close to a disc whose radius can be computed explicitely and converges to this disc at second order in~$\varepsilon$, $\varepsilon$ being the size of the roughness as well as the inverse of their wave length.

This work can be pursued in many directions: the technique used here could allow to get a better estimate of the solution (at third or further order in~$\varepsilon$). Moreover, It might be interesting to consider the same kind of study for more complex domain (not only a perturbation of a disc), or in the case of higher dimensions.
Another complexity can be introduced considering some nonlinear Bernoulli problem such that the classical $p$-Laplacian problem.

\appendix
\section{Metrics}\label{appendix}

We state here some equivalence results between distances used for interfaces.\\

Let $0<\delta<M$ be given, we consider the set of curves of $\R^2$:
$$\mathcal K=\Big\{\Gamma=\partial\Omega \mbox{ where $\Omega\subset\mathcal B(0,M)$ is a starlike domain with respect to all points in }\mathcal B(0,\delta)\Big\}.$$

Note that, for all $\Gamma\in\mathcal K$, there exists a unique parametrization $f_\Gamma:\T\rightarrow\R_+$ of $\Gamma$ in polar coordinates, which then satisfies:
$$
\Gamma=\{(f_\Gamma(\theta),\theta)\in\R^2, \theta\in \T\}.
$$
We define a first distance on the set~$\mathcal K$ based on this parametrization:
$$
d_1(\Gamma_1,\Gamma_2)=\|f_{\Gamma_2}-f_{\Gamma_1}\|_{L^{\infty}}.
$$
We recall the definition of the metric defined in  \cite{AM} {\it (which is denoted by $\Delta$ in \cite{AM})}:
$$
d_2(\Gamma_1,\Gamma_2)=\sup\{|\ln(\lambda)|, \, \lambda\in\R_+ ~ ; ~ \lambda\Gamma_1\cap\Gamma_2\neq\emptyset\}.
$$
We also recall the definition of the classical Hausdorff distance {\it (denoted here $\mathcal D_{\mathcal H}$)}:
$$
\mathcal D_{\mathcal H}(\Gamma_1,\Gamma_2)=\max\left(\sup\{d(x,\Gamma_2), x\in\Gamma_1\},\sup\{d(x,\Gamma_1), x\in\Gamma_2\}\right)
$$
where $d$ is the classical (euclidian) distance.

\begin{proposition}
The distances $d_1$, $d_2$ and $\mathcal D_{\mathcal H}$ are equivalent on $\mathcal K$.
\end{proposition}

{\bf Proof: }
Let us first remark that, for all $\Gamma_1$ and $\Gamma_2$ in $\mathcal K$, we have:
$$
d_1(\Gamma_1,\Gamma_2)\le M-\delta,
$$
$$
d_2(\Gamma_1,\Gamma_2)\le \ln(M/\delta),
$$
$$
\mathcal D_{\mathcal H}(\Gamma_1,\Gamma_2)\le M-\delta.
$$

\begin{itemize}
\item {\bf Step $1$ -} 
We prove that for all $(\Gamma_1,\Gamma_2)\in \mathcal K^2$ we have $d_1(\Gamma_1,\Gamma_2)\le \frac{M^2}{\delta} d_2(\Gamma_1,\Gamma_2)$.\\
Without loss of generality, we can assume that there exists $\theta_0\in \T$ such that:
$f_{\Gamma_2}(\theta_0)=f_{\Gamma_1}(\theta_0)+d_1(\Gamma_1,\Gamma_2)$.\\
We have $f_{\Gamma_2}(\theta_0)=\lambda f_{\Gamma_1}(\theta_0)$ with $\lambda=\displaystyle\frac{f_{\Gamma_1}(\theta_0)+d_1(\Gamma_1,\Gamma_2)}{f_{\Gamma_1}(\theta_0)}$, and then $\lambda\Gamma_1\cap\Gamma_2\neq\emptyset$ for this choice of~$\lambda$.
We then deduce:
$$
\ln\left(\frac{f_{\Gamma_1}(\theta_0)+d_1(\Gamma_1,\Gamma_2)}{f_{\Gamma_1}(\theta_0)}\right)\le d_2(\Gamma_1,\Gamma_2)
$$
and then:
$$d_1(\Gamma_1,\Gamma_2)\le f_{\Gamma_1}(\theta_0)\left(\exp(d_2(\Gamma_1,\Gamma_2))-1\right)$$
Using $f_{\Gamma_1}(\theta_0)\le M$ and $\exp(d_2(\Gamma_1,\Gamma_2))-1\le d_2(\Gamma_1,\Gamma_2)\exp(d_2(\Gamma_1,\Gamma_2))\le d_2(\Gamma_1,\Gamma_2)\times(M/\delta)$, we get the desired result.
\item {\bf Step $2$ -} 
We prove that for all $(\Gamma_1,\Gamma_2)\in \mathcal K^2$ we have $d_2(\Gamma_1,\Gamma_2)\le \frac{M}{\delta^2} \mathcal D_{\mathcal H}(\Gamma_1,\Gamma_2)$.\\
Without loss of generality, we can assume that there exists~$\theta_0\in \T$ such that
$f_{\Gamma_2}(\theta_0)=\exp(d_2(\Gamma_1,\Gamma_2))f_{\Gamma_1}(\theta_0)$.
The bounded domain $\Omega_2$ (such that $\Gamma_2=\partial\Omega_2$) being starlike with respect to all points in $\mathcal B(0,\delta)$, it contains the domain 
$\mathcal T=\cup_{m\in\mathcal B(0,\delta)}[m,f_{\Gamma_2}(\theta_0)]$, where $[m,f_{\Gamma_2}(\theta_0)]$ denotes the line delimited by these two points. Note that $\mathcal T$ is actually the union of a triangle and $\mathcal B(0,\delta)$.\\
Simple geometric arguments shows that 
$$
d(f_{\Gamma_1}(\theta_0),\partial\Gamma_2)\ge d(f_{\Gamma_1}(\theta_0),\partial\mathcal T)
$$
and 
$$
d(f_{\Gamma_1}(\theta_0),\partial\mathcal T)=\delta\frac{f_{\Gamma_2}(\theta_0)-f_{\Gamma_1}(\theta_0)}{f_{\Gamma_2}(\theta_0)}.
$$
We then get {\it (denoting $d_2=d_2(\Gamma_1,\Gamma_2)$)}:
\begin{equation*}
\begin{aligned}
d_2
& \leq \exp(d_2)-1\le\frac{\exp(d_2)-1}{\exp(d_2)}\times\exp(d_2)
\leq\frac{f_{\Gamma_2}(\theta_0)-f_{\Gamma_1}(\theta_0)}{f_{\Gamma_2}(\theta_0)}\times\exp(d_2) \\
& \leq \frac{f_{\Gamma_2}(\theta_0)-f_{\Gamma_1}(\theta_0)}{f_{\Gamma_2}(\theta_0)}\times \frac M\delta
\leq \frac{d(f_{\Gamma_1}(\theta_0),\partial\Gamma_2)}\delta\times\frac M\delta
\leq \frac M{\delta^2}\mathcal D_{\mathcal H}(\Gamma_1,\Gamma_2)
\end{aligned}
\end{equation*}
which gives the desired result.

\item {\bf Step $3$ -} 
We prove that for all $(\Gamma_1,\Gamma_2)\in \mathcal K^2$ we have $\mathcal D_{\mathcal H}(\Gamma_1,\Gamma_2)\le d_1(\Gamma_1,\Gamma_2)$.\\
Let  $x\in\Gamma_1$ such that $d(x,\Gamma_2)=\mathcal D_{\mathcal H}(\Gamma_1,\Gamma_2)$.
Noting $x=(f_{\Gamma_1}(\theta),\theta)$, $\theta\in \T$, it is clear that $|f_{\Gamma_2}(\theta)-f_{\Gamma_1}(\theta)|\ge d(x,\Gamma_2)$ which proves the desired result.

\end{itemize}

\end{document}